\documentclass[a4paper,fleqn]{cas-sc}
\usepackage{amsthm}
\usepackage{tikz}
\usetikzlibrary{arrows.meta}

\usepackage[numbers]{natbib}
\usepackage[nameinlink,noabbrev]{cleveref}

\RenewDocumentCommand\printorcid{}{}

\newtheorem{theorem}{Theorem}
\newtheorem{proposition}[theorem]{Proposition}
\newtheorem{lemma}[theorem]{Lemma}
\newtheorem{corollary}[theorem]{Corollary}
\crefname{theorem}{Theorem}{Theorems}
\Crefname{theorem}{Theorem}{Theorems}
\crefname{appendix}{appendix}{appendices}
\Crefname{appendix}{Appendix}{Appendices}

\newcommand{\icp}{\operatorname{icp}}
\newcommand{\ict}{\operatorname{ict}}
\newcommand{\oct}{\operatorname{oct}}
\newcommand{\iocp}{\operatorname{iocp}}
\newcommand{\supp}{\operatorname{supp}}

\newcommand{\ip}{\operatorname{IP}}
\newcommand{\R}{\mathbb{R}}
\newcommand{\Q}{\mathbb{Q}}
\newcommand{\Z}{\mathbb{Z}}
\newcommand{\ca}{\mathcal}
\newcommand{\size}{\mathcal{L}}
\newcommand{\bb}{\mathbb}

\def\tsc#1{\csdef{#1}{\textsc{\lowercase{#1}}\xspace}}
\tsc{WGM}
\tsc{QE}

\begin{document}
\let\WriteBookmarks\relax
\def\floatpagepagefraction{1}
\def\textpagefraction{.001}

\shorttitle{}    

\shortauthors{Yuheng Shi}  

\title [mode = title]{Integer Programming Parametrized by Imbalanced Circuit Packing}  



%

\author[1]{Yuheng Shi}

\cormark[1]


\ead{s22yshi@uni-bonn.de}



\affiliation[1]{organization={University of Bonn},
            addressline={Regina-Pacis-Weg 3},
            city={Bonn},
            postcode={53113},
            country={Germany}}






\cortext[1]{Corresponding author}



\begin{abstract}
We study the standard form Integer Programming problem (IP) $\min\{ c^t x:Ax=b,x\in \mathbb{Z}^n_{\geq 0}\}$ parametrized by the \emph{imbalanced circuit packing number}, $\icp(A)$, the maximum number of support-minimal nonzero vectors in $\ker(A)$ such that these vectors have pairwise-disjoint support and cannot be scaled to $\{-1,0,1\}$-vectors by nonzero scalars. This packing number is $0$ if and only if $A$ can be transformed to a totally unimodular matrix by elementary row operations, in which case (IP) is solvable in polynomial time. We show that if $A$ represents a graphic matroid, (IP) is fixed-parameter tractable when parametrized by $\icp(A)$.The proof combines an  Erd\H{o}s--P\'osa property for edge-colored non-monochromatic cycles with a reduction to mixed-integer programming using few integer variables. In contrast, by a reduction from the parity potential problem, we show NP-completeness of feasibility of (IP) when $A$ represents a cographic matroid, $\ker(A)$ is $2$-regular and $\icp(A)\leq 1$.

\end{abstract}



\begin{keywords}
Integer programming \sep Circuit imbalance measure \sep Fixed parameter tractability \sep Graphic matroids \sep Cographic matroids \sep Erd\H{o}s--P\'osa property
\end{keywords}

\maketitle

\section{Introduction}\label{intro}

We consider the standard-form integer program
\begin{equation}\label{eq:ip}
 (\ip)\qquad \min\{c^\top x:Ax=b,\ x\in\Z^n_{\geq0}\},
\end{equation}
where $A\in\Z^{m\times n}$, $b\in\Z^m$, and $c\in\Z^n$. It is well-known that $\eqref{eq:ip}$ is solvable in polynomial time if $A$ is totally unimodular (TU), i.e. all subdeterminants of $A$ are in $\{-1,0,1\}$. To generalize this result, it is natural to ask whether \eqref{eq:ip} is still solvable in polynomial time when $A$ is totally $\Delta$-modular, i.e. all subdeterminants of $A$ are bounded in absolute value by $\Delta$ for some fixed $\Delta\in\Z_{\geq 1}$. It is known that there is a strongly polynomial algorithm when $\Delta = 2$ (see \cite{ArtmannWeismantelZenklusen2017}). It is an open problem whether there is a polynomial time algorithm when $\Delta \geq 3$, and various special cases have been studied (see  \cite{FioriniJoretWeltgeYuditsky2025,Kober2026,AprileFioriniJoretKoberSewerynWeltgeYuditsky2026,NageleNobelSantiagoZenklusen2025}).

In this paper, we study another way of generalizing TU-optimization from the \emph{circuit perspective} (see \cite{EkbataniNaturaVegh2022} for relevant discussions). Specifically, for any subspace $W\subseteq \R^n$, $v\in W$ is an \emph{elementary vector} if $v$ is nonzero and  has minimal support among all nonzero vectors $w\in W$. That is, $v\neq 0$ and if $w\in W-\{0\}$ satisfies $\supp(w)\subseteq \supp (v)$, then $\supp(w)=\supp(v)$, where $\supp(v):=\{i\in [n]: v_{i}\neq 0\}$. $C\subseteq [n]$ is a \emph{circuit} of $W$ if it is the support of an elementary vector of $W$. We denote the set of circuits of $W$ by $\mathcal{C}_W$. Note that if $W=\ker A$, then circuits of $W$ are exactly circuits of the vector matroid $M(A)$. Furthermore, we say that an elementary vector $v\in W$ is balanced if $|v_i|=|v_j|$ for all $i,j\in\supp(v)$. and imbalanced otherwise. We say $C\in\ca{C}_W$ is balanced if it is the support of a balanced elementary vector, and imbalanced otherwise. Note that if $v_1,v_2\in W$ satisfy $\supp(v_1)=\supp(v_2)\in\ca{C}_W$, then $v_1=\lambda v_2$ for some $\lambda\in\R\backslash \{0\}$, so the definitions of a balanced circuit and an imbalanced circuit make sense.

There is a relationship between TU and non-existence of imbalanced circuits, made precise by the following two statements (see Theorem 3.4 of \cite{EkbataniNaturaVegh2022} for a proof): first, if $A$ is TU, then $\ker(A)$ has no imbalanced circuits; second, if $A\in\R^{m\times n}$ is in basis form (i.e. containing an $m\times m$ identity matrix) and $\ker(A)$ has no imbalanced circuits, then $A$ is TU. Because we can use Gaussian elimination and remove zero rows (if any) to transform any $A\in\R^{m\times n}$ into a basis form (and apply the same row operations to $b$), and because these operations do not change $\ker(A)$, we can solve \eqref{eq:ip} in polynomial time if $A$ has no imbalanced circuits.

This leads to the following natural definition and question: for any $A\in\R^{m\times n}$, let $\icp(A)$, the \emph{imbalanced circuit packing number} of $A$, denote the maximum number of pairwise disjoint imbalanced circuits of $A$. Is $\ip$ still solvable in polynomial time when $\icp(A)\leq k$ for some fixed $k\in\Z_{\geq 1}$?

Before answering this question, let us motivate the definition $\icp(A)$ in three different ways.

First, as shown by the previous discussion, \eqref{eq:ip} is solvable in polynomial time if $\icp(A)=0$, which is equivalent to $A$ being TU up to elementary row operations.

Second, a transversal set for imbalanced circuits gives an efficient algorithm for \eqref{eq:ip}. Specifically, for $A\in\R^{m\times n}$, we say $S\subseteq[n]$ is an \emph{imbalanced circuit transversal} of $A$ if $S\cap C\neq\emptyset$ for any imbalanced circuit $C$ of $A$, and let $\ict(A)$, the \emph{imbalanced circuit transversal number}, be the smallest size of an imbalanced circuit transversal of $A$. As long as an imbalanced circuit transversal $S\subseteq [n]$ is given, we can reduce \eqref{eq:ip} to a mixed-integer linear program with exactly $|S|$ integer variables, which can be solved in polynomial time if $|S|$ is fixed, see \Cref{prop:milp}. We obviously have $\icp(A)\leq \ict(A)$. It turns out that when $A$ represents graphic matroids, $\ict(A)\leq c\cdot\icp(A)(\icp(A)+1)\log(\icp(A)+1)$ for some $c>0$ (see \Cref{EP3}), so an Erd\H{o}s--P\'osa property holds and \eqref{eq:ip} can be efficiently solved (\Cref{thm:graphic}). This does not hold already for cographic matrices (see \Cref{EP4}).

Third, there is a connection to the minimum vertex cover problem (or equivalently the maximum stable set problem). For an undirected graph $G=(V,E)$, let $A_G$ be the unsigned vertex-edge incidence matrix of $G$, and put $M_G=[A_G^t|-I_{m}]$ where $m=|E(G)|$. For $c\in\Q_{\geq 0}^{V(G)}$, the \emph{minimum weight vertex cover problem} can be expressed as
\[
 \min\{c^\top y:M_G(y,s)^\top=\mathbf{1}_{m},\ (y,s)\in\Z_{\geq0}^{V\cup E}\}.
\]

We prove in \Cref{vc} that
\begin{equation}\label{eq:graph-identity}
 \icp(M_G)=\iocp(G),\qquad \ker(M_G) \text{ is 2-regular},
\end{equation}
where $\iocp(G)$, standing for \emph{induced odd cycle packing} number, is the maximum number of vertex-disjoint odd cycles with no edge joining distinct cycles. A subspace $W\subseteq \R^n$ is called $2$-regular if all elementary vectors can be scaled to a $\{-2,-1,0,1,2\}$-vector by a nonzero scalar (see \cite{Lee1990} for more on $k$-regularity). The minimum vertex cover problem on graphs with a bounded number of vertex-disjoint odd cycles is polynomial time solvable, see \cite[Theorem 2]{FioriniJoretWeltgeYuditsky2025}. But the complexity of minimum vertex cover on graphs $G$ with $\iocp(G)\leq1$ remains open and solving this problem would settle the long-standing open problems whether the maximum clique problem is polynomial time solvable for disk graphs and unit ball graphs; see introduction of \cite{BonamyEtAl2024}. Our hardness result \Cref{thm:cographic} shows that optimizing the more general class of integer programs satisfying the two bounds in \eqref{eq:graph-identity}, i.e. integer programs of form \eqref{eq:ip} where $\icp(A)\leq 1$ and $\ker(A)$ is 2-regular, is NP-hard.

Our main results, \Cref{thm:graphic} and \Cref{thm:cographic}, describe tractability and hardness of \eqref{eq:ip} parametrized by $\icp(A)$ in terms of the matroid structure of $A$.

\begin{theorem}\label{thm:graphic}
Suppose that $A$ represents a graphic matroid and let $k\in\Z_{\geq0}$. In time
$$
 2^{O(k^2\log^2(k+1))}\size^{O(1)},
$$
one can either find $k+1$ pairwise disjoint imbalanced circuits or solve \eqref{eq:ip}, where $\size$ is the binary encoding length of the instance \eqref{eq:ip}. In particular, integer programming is fixed-parameter tractable in $\icp(A)$ when $A$ represents a graphic matroid.
\end{theorem}


We note that deciding whether $\{x\in\Z^n_{\geq0}:Ax=b\}$ is nonempty is NP-complete when $A$ represents a graphic matroid: indeed, consider the class of matrices $A$ where $A$ is the vertex-edge incidence matrix of a directed graph $G$ with some columns $S\subseteq [n]$ scaled by $2$. Then $A$ represents the cycle matroid of the underlying undirected graph of $G$, and $\{x\in\Z^n_{\geq0}:Ax=b\}$ is nonempty if and only if there is an uncapacitated $b$-flow $f:E(G)\to \mathbb{Z}_{\geq 0}$ on $G$ such that $f(e)$ is even for $e\in S$. This is the \emph{parity network flow problem}, known to be NP-complete (see Section 3 of \cite{Conforti2009Network}). Consider a more general problem, where given divisibility constraints $\{d_e\}_{e\in E(G)}$, $d_e\in \bb{Z}_{>0}$, we ask for an integral $b$-flow $f$ of minimum cost such that $f_e$ is required to be in $d_e\mathbb{Z}$ for all $e\in E(G)$. \Cref{thm:graphic} states that this problem is solvable in polynomial time if we fix the maximum number of edge-disjoint (undirected) cycles $C$ where $d_e\neq d_{e'}$ for some $e,e'\in E(C)$.

\begin{theorem}\label{thm:cographic}
Deciding whether $\{x\in\Z^n_{\geq0}:Ax=b\}$ is nonempty is NP-complete even for integral $A,b$ such that $A$ represents a cographic matroid, $\ker(A)$ is $2$-regular, and $\icp(A)\leq 1$.
\end{theorem}

The proof follows from NP-completeness of the natural dual problem to the parity network flow problem, which we call \emph{parity potential problem}, see \Cref{HardnessParityPotential}.

\section{FPT parametrized by $\icp(A)$ when $A$ is graphic}
For $T\subseteq[n]$, write $A_T$ for the columns indexed by $T$. Elementary vectors of $\ker(A_T)$, extended by zeros, are elementary vectors of $\ker(A)$. Consequently, if $S$ is an imbalanced circuit transversal and $T=[n]\setminus S$, every circuit of $A_T$ is balanced. After Gaussian elimination, nonzero rows of $A_T$ can be put in basis form, that is, containing an identity matrix of size $\mathrm{rank}(A_T)$. The circuit characterization of total unimodularity then implies that the resulting $A_T$ is TU \cite[Theorem~3.4]{EkbataniNaturaVegh2022}.

To prove \Cref{thm:graphic}, we need the fact that $\eqref{eq:ip}$ is efficiently solvable when an imbalanced circuit transversal of bounded size is given:

\begin{proposition}\label{prop:milp}
Given an instance of \eqref{eq:ip} and an imbalanced circuit transversal $S$ of size $k$, one can, in polynomial time, either conclude infeasibility or construct an equivalent mixed-integer linear program with $k$ integer variables. Thus \eqref{eq:ip} can be solved in time $2^{O(k\log(k+1))}\size^{O(1)}$, where $\size$ is the binary encoding length of \eqref{eq:ip}.
\end{proposition}

The reduction to mixed integer linear programming uses Smith normal form and solves linear congruence equations. Then the mixed-integer linear program is solved by an algorithm of Kannan (\cite{Kannan1987Minkowski}). The complete proof will be deferred to \Cref{Prop3Proof}.

The second ingredient of the proof of $\cref{thm:graphic}$ is an Erd\H{o}s--P\'osa property for $\icp(A)$ and $\ict(A)$ when $A$ represents a graphic matroid. For an undirected graph $G$, $S\subseteq V(G)$, and $T\subseteq E(G)$, an ``$S$-cycle'' (``$T$-cycle'') means a cycle $C$ in $G$ such that $V(C)\cap S\neq\emptyset$ ($E(C)\cap T\neq\emptyset$). Our starting point is the following result:

\begin{lemma}\cite[Theorem~2]{PontecorviWollan2012}\label{EP0}
    There exists a constant $c>0$ such that for all graphs $G$, $S\subseteq V(G)$, and positive integers $k$, either
there exist $k$ vertex disjoint $S$-cycles, or there exists at most $ck \log k$ vertices hitting every $S$-cycle. There
exists an $O(nm)$-time algorithm which given in input $G$, $S$ and $k$, returns either $k$ vertex disjoint $S$-cycles or
an $S$-cycle vertex hitting set of size at most $ck \log k$.
\end{lemma}

The next result is an ``edge version'' of \Cref{EP0}:

\begin{lemma}\label{EP1}
    There exists a constant $c>0$ such that, for each positive integer $k$, undirected graph $G$ and $T\subseteq E(G)$, either $G$ has $k$ edge-disjoint $T$-cycles, or there exists at most $ck\log(k)$ edges hitting every $T$-cycle. There is an algorithm returning either $k$ edge-disjoint $T$-cycles or at most $ck\log(k)$ edges hitting every $T$-cycle, with running time $\mathcal{O}(m^3n)$ where $m=|E(G)|$ and $n=|V(G)|$.
\end{lemma}
\begin{proof}
    Consider edges in $T$ as colored edges. Consider the following reduction: replace each edge of $G$ by a path of length $3$, color the middle edge if the original edge is colored, then take the line graph of this graph and color a vertex of this line graph if its corresponding edge is colored. Call this new graph $H$, with colored vertices $S\subseteq V(H)$. Then the maximum number of edge-disjoint $T$-cycles in $G$ is equal to the maximum number vertex-disjoint $S$-cycles in $H$. Note that to get a $T$-cycle in $G$ from an $S$-cycle in $H$, we might need to skip some vertices of the $S$-cycle which correspond to a star in $G$.
    
    If $G$ does not have $k$ edge-disjoint $T$-cycles, then $H$ does not have $k$ vertex-disjoint $S$-cycles, then \Cref{EP0} implies $H$ has at most $ck\log(k)$ vertices intersecting each $S$-cycle for some constant $c$, which correspond to at most $ck\log(k)$ edges of $G$ intersecting each $T$-cycle. The algorithmic result follows because \Cref{EP0} states an algorithm. The running time follows because $|V(H)|=3|E(G)|$ and $|E(H)|=\mathcal{O}(|E(G)|^2\cdot|V(G)|)$.
\end{proof}

We use \Cref{EP1} to prove an Erd\H{o}s--P\'osa property for non-monochromatic cycles in the next result. Given a graph $G$ and a function $\chi:E(G)\to \Z_{> 0}$ which we also refer to as one (improper) coloring of edges, we say a cycle $C$ in $G$ is non-monochromatic if there exist $e,f\in E(C)$ such that $\chi(e)\neq\chi(f)$. Non-monochromatic cycles have been studied in different contexts, for example in \cite{IwataYamaguchiYoshida2018} and \cite{Wahlstrom2017}, where it is shown that there is a FPT algorithm to find at most $k$ vertices $V$ such that $G\backslash V$ has no non-monochromatic cycles (or decides there is nonexx), where the parameter is $k\in\bb{N}$. The following result investigates the relationship between (edge) packing and hitting edge set for non-monochromatic cycles, which has not been studied before to the best of our knowledge.
\begin{lemma}\label{EP2}
There exists a constant $c>0$ such that, for each positive integer $k$, undirected graph $G$ and a coloring of edges $\chi:E(G)\to\bb{Z}_{> 0}$, either $G$ has $k$ edge-disjoint non-monochromatic cycles, or there exists at most $c(k-1)k\log(k)$ edges hitting every non-monochromatic cycle. There is an algorithm returning either $k$ edge-disjoint non-monochromatic cycles or at most $c(k-1)k\log(k)$ edges hitting every non-monochromatic cycle, with running time polynomial in $m=|E(G)|$ and $n=|V(G)|$.
\end{lemma}
\begin{proof}
    We replace each vertex $v$ of $V(G)$ by a star with $|\chi(\delta_G(v))|$ edges, and connect original edges in $E(G)$ to star leaves of the corresponding color. Formally, construct a new graph $H$ in the following way: Let $V(H)=\{v_0 : v\in V(G)\}\cup\{v_i:v\in V(G), i\in \chi(\delta_{G}(v))\}$ and let $E(H)=\{\{v_0,v_i\}:v\in V(G), i\in \chi(\delta_G(v))\}\cup \{\{u_{\chi(e)},v_{\chi(e)}\}:e=\{u,v\}\in E(G)\}$ (parallel edges of $G$ of the same color are still parallel edges in $H$). 

    Let $T=\cup_{v\in V(G)}\delta_{H}(v_0)$ be the star edges. We note that if $C\subseteq E(H)$ is an $T$-cycle in $H$, then after viewing $C\backslash T$ as an edge set in $G$, each vertex in $(V(G),C\backslash T)$ has even degree. So $C\backslash T$ can be decomposed into edge-disjoint union of cycles in $G$. We claim that one of these cycles must be non-monochromatic. Indeed, because $C$ is an $T$-cycle in $H$, $C$ contains two edges of the form $\{v_0,v_i\},\{v_i,u_i\}\in E(H)$. This implies $\{v,u\}\in E(G)$ and $\chi(\{v,u\})=i$. Furthermore, $\{v,u\}$ is the only edge of $C\backslash T\subseteq E(G)$ which is adjacent to $v$ and has color $i$. Now if $C\backslash T$ can be decomposed into edge-disjoint union of monochromatic cycles in $G$, then the number of edges of $C\backslash T\subseteq E(G)$ which are adjacent to $v$ and have color $i$ must be even, contradiction. As a result, the maximum number of edge-disjoint $T$-cycles in $H$ is at most the maximum number of edge-disjoint non-monochromatic cycles in $G$.

    Now run the algorithm of \Cref{EP1} on $H$ to see whether there are $k$ edge-disjoint $T$-cycles or at most $ck\log k$ edges hitting every $T$-cycle for some constant $c>0$. In the first case, the argument of the previous paragraph gives us $k$ edge-disjoint non-monochromatic cycles in $G$. In the second case, call the $ck\log k$ edges $R$. We want to put $R$ back to $E(G)$, so we need to deal with the edges in $R\cap T$. Now take $e=\{u_0,u_i\}\in R\cap T$, and define $H_e$ by removing $e$ from $H$ and replacing each edge in $T\backslash\{e\}$ by $k$ parallel edges. To ensure polynomial running time, we can assume $k< E(G)$, because otherwise we can choose $E(G)$ to be an edge set hitting every non-monochromatic cycle. Compute the minimum $u_0$-$u_i$ cut in $H_e$, where each edge has capacity $1$. If this cut has value at least $k$, then by Menger's theorem we can find $k$ edge-disjoint $u_0$-$u_i$ paths in $H_e$ by a flow computation. Each of these paths becomes an $T$-cycle after adding $e$ back, which further contains a non-monochromatic cycle in $G$ after removing edges in $T$. So we find $k$ edge-disjoint non-monochromatic cycles in $G$. 
    
    Otherwise, we have a $u_0$-$u_i$ cut in $H_e$ of value less than $k$, which implies that all edges in this cut come from $G$, as we have replaced each edge in $T\backslash\{e\}$ with $k$ parallel edges. For each $e\in T\cap R$, remove $e$ from $R$ add these less than $k$ cut edges to $R$, to get a new set $R'\subseteq E(H)$ such that $|R'|\leq c(k-1)k\log k$ and $R'\cap T= \emptyset$. By construction, $R'$ hits every $T$-cycle of $H$. Now we can view $R'$ as a subset of $E(G)$, and $R'$ hits every non-monochromatic cycle of $G$ because such a cycle can be lifted to a $T$-cycle in $H$.

    The claimed running time also follows from the construction of $H$ and from the fact that we can compute minimum cut and flow decomposition in polynomial time.
\end{proof}

We are now ready to prove \Cref{thm:graphic}. For an undirected graph $G$, $M(G)$ denotes its cycle matroid. For a real matrix $A$, $M(A)$ denotes its vector matroid.

\emph{Proof of \Cref{thm:graphic}.}\label{proofof1} By assumption, $A$ represents a graphic matroid. Because we have a polynomial time independence oracle for this matroid by Gaussian elimination, we can find an undirected graph $G$ such that $M(G)=M(A)$ by the algorithm in \cite{Seymour1981Graphic}. Let $A_{\overrightarrow{G}}$ be the vertex-edge incidence matrix of some orientation $\overrightarrow{G}$ of $G$, then $M(A)=M(A_{\overrightarrow{G}})$ (\cite[Lemma 5.1.3]{Oxley2011Matroid}). Because binary matroids are uniquely representable over any field (\cite[Proposition 6.6.5]{Oxley2011Matroid}), there exists a rational diagonal nonsingular matrix $D\in\bb{Q}^{n\times n}$ such that $\ker(A)=D\ker(A_{\overrightarrow{G}})=\ker(A_{\overrightarrow{G}}D^{-1})$. Such $D$ can be computed in polynomial time, for example by \cite[Lemma 6.6.4]{Oxley2011Matroid} and the proof of \cite[Theorem 6.4.7]{Oxley2011Matroid}.

Let $\chi:E(G)\to\bb{Z}_{>0}$ be an edge coloring such that $\chi(e)=\chi(f)$ if and only if $|D_e|=|D_f|$ for all $e,f\in E(G)$, then imbalanced circuits of $A_{\overrightarrow{G}}D^{-1}$ correspond to non-monochromatic cycles of $G$ (To see this, first note that circuits of $A_{\overrightarrow{G}}$ are all balanced, and they correspond to cycles of $G$). Applying the algorithm of \Cref{EP2} (with input $k+1$ and $G$), we find in polynomial time either $k+1$ edge-disjoint non-monochromatic cycles of $G$, which correspond to $k+1$ pairwise disjoint imbalanced circuits of $A_{\overrightarrow{G}}D^{-1}$, which are  $k+1$ pairwise disjoint imbalanced circuits of $A$ because $\ker(A)=\ker(A_{\overrightarrow{G}}D^{-1})$, or find $\ca{O}(k^2\log (k+1))$ edges of $G$ hitting every non-monochromatic cycle. Because $\ker(A)=\ker(A_{\overrightarrow{G}}D^{-1})$, these edges correspond to an imbalanced circuit transversal of $\ker(A)$ of size $\ca{O}(k^2\log(k+1))$. By \Cref{prop:milp}, the IP $\min \{c^tx:Ax = b, x\in\bb{Z}^n_{\geq 0}\}$ can be solved in time $ 2^{\ca{O}(k^2\log^2 (k+1))}\cdot \size^{\ca{O}(1)}$. \qed

\Cref{EP2} also gives rise to the following Erd\H{o}s--P\'osa property:

\begin{corollary}\label{EP3}
    There exists a constant $c>0$ such that $\ict(A)\leq c
    \cdot\icp(A)(\icp(A)+1)\log(\icp(A)+1)$ for all $A\in\bb{Z}^{m\times n }$ representing a graphic matroid.  
\end{corollary}

\section{NP-completeness when $A$ is cographic, $\icp(A)\leq 1$, and $\ker(A)$ is $2$-regular}
To prove \Cref{thm:cographic}, we first prove NP-completeness of a natural ``dual'' problem to the parity network flow problem considered in Section 3 of \cite{Conforti2009Network}, which we call the \emph{parity potential problem}: Given a directed graph $G$, cost function $c:E(G)\to\mathbb{Z}$ and $S\subseteq E(G)$, does there exist $\pi:V(G)\to \mathbb{Z}$ such that $c^\pi(u,v):= c(u,v)+\pi(u)-\pi(v)\geq 0$ for all $e=(u,v)\in E(G)$ and $c^\pi(e)$ is even for all $e\in S$? Note that this problem is solvable in polynomial time if there is no parity constraint, because it can then be formulated as a feasibility integer programming problem, where the constraint matrix is the transpose of the vertex-edge incidence matrix of $G$, which is TU.

\begin{theorem}\label{HardnessParityPotential}
    The parity potential problem is NP-Complete, even if $c$ takes value in $\{-1, 1, 2,3\}$ and the underlying undirected graph of $G$ is complete.
\end{theorem}
\begin{proof}
    The parity potential problem is in NP, because it can be formulated as an integer linear program, and feasibility of integer linear programs is in NP.

    Consider the not-all-equal 3-satisfiability problem, abbreviated as NAE3SAT: Given a set $X$ of variables and $\mathcal{C}$ of clauses where each clause has exactly $3$ literals, does there exist a truth assignment of $X$ such that each clause contains at least one true literal and at least one false literal? This problem is known to be NP-complete (\cite{Schaefer1978Satisfiability}).

    Given an instance $(X,\mathcal{C})$ of NAE3SAT, we create the following instance $(G,c,S)$ of the parity potential problem: Add variable vertices $v_x\in V(G)$ for each $x\in X$. Add clause vertices $\{a_C, b_C, c_C\}\subseteq V(G)$ for each $C=\{a,b,c\}\in\mathcal{C}$. These are all vertices of $G$. 
    
    For each clause $C=\{a,b,c\}$, add edges $(a_C, b_C), (b_C, c_C), (c_C, a_C)$ with costs $c(a_C,b_C)=1, c(b_C,c_C)=1, c(c_C,a_C)=-1$. These edges are not in $S$. 

    Add edges from variable vertices to clause vertices according to occurrence in the following way. For all $x\in X$ and $C=\{a,b,c\}\in \mathcal{C}$: 
    
    add an edge $(v_x, a_C)$ with $c(v_x, a_C)=2$ if $a=x$; 
    
    add an edge $(v_x, a_C)$ with $c(v_x, a_C)=3$ if $a=\overline{x}$; 
    
    add an edge $(v_x, b_C)$ with $c(v_x, b_C)=2$ if $b=x$; 
    
    add an edge $(v_x, b_C)$ with $c(v_x, b_C)=3$ if $b=\overline{x}$; 
    
    add an edge $(v_x, c_C)$ with $c(v_x, c_C)=2$ if $c=x$; 
    
    add an edge $(v_x, c_C)$ with $c(v_x, c_C)=3$ if $c=\overline{x}$. 
    
    All edges added in these ways are in $S$.
    Finally, make the underlying undirected graph complete by adding edges of cost $3$, oriented arbitrarily. These edges are not in $S$. This finishes the construction of the parity potential problem instance $(G,c,S)$. Note that $c$ indeed only takes values in $\{-1, 1, 2, 3\}$.

    See \Cref{fig:nae3sat-parity-potential} for an illustration of the construction.

    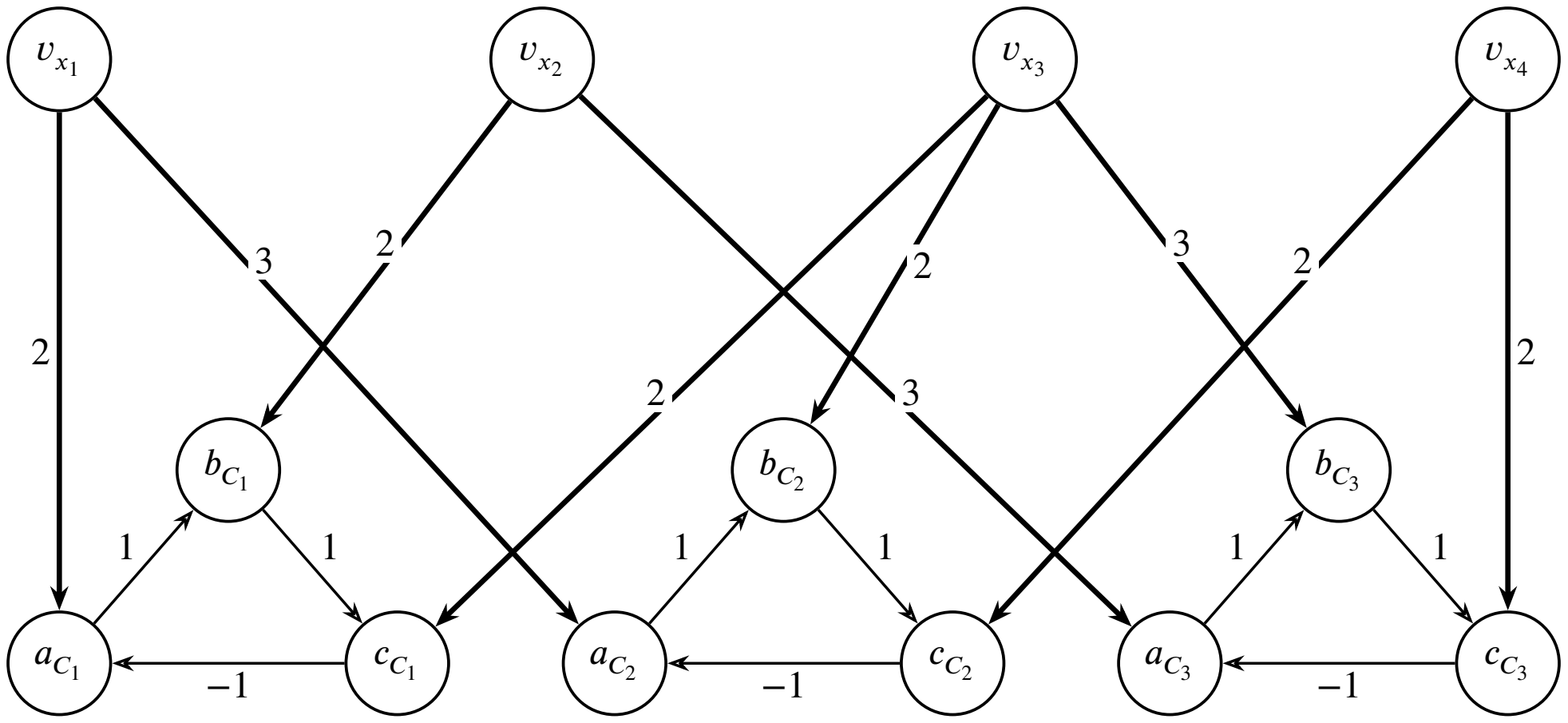
\begin{figure}[pos=!t]
\centering
\begin{tikzpicture}[
    x=1cm, y=1cm,
    vertex/.style={circle, draw, fill=white, line width=.7pt,
                   minimum size=8.5mm, inner sep=1pt, font=\small},
    arc/.style={line width=.65pt,
                -{Stealth[open,length=1.8mm,width=1.3mm]}},
    special/.style={line width=1.25pt,
                    -{Stealth[length=2mm,width=1.5mm]}},
    cost/.style={font=\small, fill=white, inner sep=1.5pt}
]
\foreach \i/\x in {1/0,2/4,3/8,4/12}
    \node[vertex] (v\i) at (\x,5) {$v_{x_\i}$};

\foreach \i/\x in {1/0,2/4.6,3/9.2} {
    \node[vertex] (a\i) at (\x,0)       {$a_{C_\i}$};
    \node[vertex] (b\i) at (\x+1.4,1.6) {$b_{C_\i}$};
    \node[vertex] (c\i) at (\x+2.8,0)   {$c_{C_\i}$};
}

\draw[special] (v1) -- node[cost,pos=.48,left]  {$2$} (a1);
\draw[special] (v2) -- node[cost,pos=.50,above] {$2$} (b1);
\draw[special] (v3) -- node[cost,pos=.60,above] {$2$} (c1);

\draw[special] (v1) -- node[cost,pos=.35,above] {$3$} (a2);
\draw[special] (v3) -- node[cost,pos=.50,right] {$2$} (b2);
\draw[special] (v4) -- node[cost,pos=.35,above] {$2$} (c2);

\draw[special] (v2) -- node[cost,pos=.60,above] {$3$} (a3);
\draw[special] (v3) -- node[cost,pos=.50,above] {$3$} (b3);
\draw[special] (v4) -- node[cost,pos=.48,right] {$2$} (c3);

\foreach \i in {1,2,3} {
    \draw[arc] (a\i) -- node[cost,midway,above left]  {$1$} (b\i);
    \draw[arc] (b\i) -- node[cost,midway,above right] {$1$} (c\i);
    \draw[arc] (c\i) -- node[cost,midway,below]       {$-1$} (a\i);
}
\end{tikzpicture}
\caption{The construction for the NAE3SAT instance with
$X=\{x_1,x_2,x_3,x_4\}$ and $\mathcal{C}=\{C_1,C_2,C_3\}$, where
$C_1=\{x_1,x_2,x_3\}$,
$C_2=\{\overline{x_1},x_3,x_4\}$, and
$C_3=\{\overline{x_2},\overline{x_3},x_4\}$.
Each clause must contain both a true and a false literal for this instance to be a yes-instance.
Thick arcs belong to $S$; edge labels give costs.
The edges added to complete the graph are omitted.}
\label{fig:nae3sat-parity-potential}
\end{figure}
    Note that if there exists $\pi:\{a_C, b_C, c_C\}\to \mathbb{Z}$ such that $c^{\pi}$ is nonnegative on $(a_C, b_C), (b_C, c_C)$ and $(c_C, a_C)$, then $(\pi(a_C),\pi(b_C), \pi(c_C))\not\equiv (0,0,0),(1,1,1) \mod 2$. Indeed, if $(\pi(a_C),\pi(b_C), \pi(c_C))\equiv (0,0,0)\mod 2$ or $(\pi(a_C),\pi(b_C), \pi(c_C))\equiv (1,1,1)\mod 2$, then $c^{\pi}(a_C,b_C)$, $c^{\pi}(b_C,c_C)$, and $c^{\pi}(c_C,a_C)$ are odd, so nonnegativity of $c^\pi$ implies $c^\pi$ is at least $1$ on these three edges, so $1=c(a_C,b_C)+c(b_C,c_C)+c(c_C,a_C)=c^\pi(a_C,b_C)+c^\pi(b_C,c_C)+c^\pi(c_C,a_C)\geq 3$, contradiction.
    
    Conversely, for each parity pattern other than $(0,0,0) \mod 2$ and $(1,1,1) \mod 2$, there exists $\pi:\{a_C,b_C,c_C\}\to \{-1, 0, 1, 2\}$ such that $(\pi(a_C),\pi(b_C), \pi(c_C))$ has this parity pattern and $c^\pi$ is nonnegative on $(a_C, b_C), (b_C, c_C)$ and $(c_C, a_C)$. The following is a complete list of feasible $\pi$-values according to the $6$ different parity patterns: 
    
    for $(\pi(a_C), \pi(b_C), \pi(c_C))\equiv (1,0,0) \mod 2$, let $(\pi(a_C), \pi(b_C), \pi(c_C))=(-1,0,0)$; 
    
    for $(\pi(a_C), \pi(b_C), \pi(c_C))\equiv (0,1,0) \mod 2$, let $(\pi(a_C), \pi(b_C), \pi(c_C))=(0,1,2)$; 
    
    for $(\pi(a_C), \pi(b_C), \pi(c_C))\equiv (0,0,1) \mod 2$, let $(\pi(a_C), \pi(b_C), \pi(c_C))=(0,0,1)$; 
    
    for $(\pi(a_C), \pi(b_C), \pi(c_C))\equiv (1,1,0) \mod 2$, let $(\pi(a_C), \pi(b_C), \pi(c_C))=(-1,-1,0)$; 
    
    for $(\pi(a_C), \pi(b_C), \pi(c_C))\equiv (1,0,1) \mod 2$, let $(\pi(a_C), \pi(b_C), \pi(c_C))=(-1,0,1) $; 
    
    for $(\pi(a_C), \pi(b_C), \pi(c_C))\equiv (0,1,1) \mod 2$, let $(\pi(a_C), \pi(b_C), \pi(c_C))=(0,1,1)$.
    
    Now if the NAE3SAT instance $(X, \mathcal{C})$ is a yes-instance with a valid truth assignment $T:X\to\{\mathrm{true}, \mathrm{false}\}$, define $\pi:\{v_x: x\in X\}\to \{0,1\}$ by $\pi(v_x)=1$ if $T(x)=\mathrm{true}$ and $\pi(v_x)=0$ if $T(x)=\mathrm{false}$. For each clause $C=\{a,b,c\}$, we know $(T(a),T(b), T(c))\neq (\mathrm{true},\mathrm{true},\mathrm{true})$ and $(T(a),T(b), T(c))\neq (\mathrm{false},\mathrm{false},\mathrm{false})$. Thus we can find $\pi:\{a_C,b_C,c_C\}\to \{-1,0,1,2\}$ such that $c^{\pi}$ is nonnegative on $(a_C,b_C), (b_C, c_C), (c_C, a_C)$ and $(\pi(a_C), \pi(b_C), \pi(c_C))\equiv (T(a), T(b), T(c)) \mod 2$ for all clauses $C=\{a,b,c\}$, where we interpret ``true'' as $1$ and interpret ``false'' as $0$. Note that $c^\pi(e)\geq 0$ for all $e\in E(G)$ because $\pi(v_x)\in\{0, 1\}$ for all $x\in X$ and $\pi(v)\in\{-1, 0, 1,2\}$ for all $v\in V(G)$. It is not hard to see that $c^\pi(e)\equiv 0 \mod 2$ for all $e\in S$.

    Conversely, if $(G,c,S)$ is a yes-instance, i.e. there is $\pi:V(G)\to \mathbb{Z}$ such that $c^\pi(e)\geq 0$ for all $e\in E(G)$ and $c^\pi(e)\equiv 0 \mod 2$ for all $e\in S$, then define $T:X\to \{\mathrm{true}, \mathrm{false}\}$ by $T(x)=true$ if $\pi(v_x)\equiv 1\mod 2$ and $T(x)=false$ if $\pi(v_x)\equiv 0\mod 2$. Because we know that for each clause $C=\{a,b,c\}$, $(\pi(a_C),\pi(b_C),\pi(c_C)) \not\equiv (0,0,0)\mod 2$ and $(\pi(a_C),\pi(b_C),\pi(c_C)) \not\equiv (1,1,1)\mod 2$, $T$ has to be a valid truth assignment for the NAE3SAT instance $(X, \mathcal{C})$.
\end{proof}

We need the following terminology for the proof of \Cref{thm:cographic}: for a matrix $A$, $M(A)$ will denote the vector matroid of $A$ and $r.s.(A)$ will denote the row space of $A$. For an undirected graph $G$, $M(G)$ will denote the cycle matroid of $G$. For a matroid $M$, $M^*$ will denote the dual matroid of $M$. See \cite{Oxley2011Matroid} for discussion of these concepts.

\emph{Proof of \Cref{thm:cographic}.} The problem of deciding whether $\{x\in\mathbb{Z}^n: Ax=b, x\geq 0\}$ is nonempty is in NP because integer linear programming is in NP. We prove it is NP-complete by a reduction from \Cref{HardnessParityPotential}.

Let $(G, c, S)$ be an instance of the parity potential problem, where the underlying undirected graph of $G$, which we denote by $\tilde{G}$, is complete. Write $n:=|V(G)|$ and $m:=|E(G)|$. Let $A$ be the incidence matrix of $G$, then $M(A)=M(\tilde{G})$ (\cite[Lemma 5.1.3]{Oxley2011Matroid}). By restricting to a maximal set of linearly independent rows of $A$ (which is in fact any $n-1$ rows of $A$), then picking a basis and multiplying the inverse of this basis on the left, we get a real matrix $A'$ such that $\ker(A')=\ker(A)$, $r.s.(A')=r.s.(A)$, $M(A')=M(A)$, and $A'$ is in basis form. Because $\ker(A')=\ker(A)$ and $A$ is totally unimodular, we know by \cite[Theorem 3.4]{EkbataniNaturaVegh2022} that $A'$ is totally unimodular. Write $A'=[I_{n-1} | N]$, and let $B=[-N^{t} |I_{m-(n-1)}]$. By \cite[Theorem 2.2.8]{Oxley2011Matroid}, $M(B)=M(A')^*=M(A)^*=M(\tilde{G})^*$, so $B$ represents a cographic matroid. Furthermore, $\ker(B)=r.s.(B)^\perp = r.s.(A') = r.s.(A)$ (by \cite[Proposition 2.2.23]{Oxley2011Matroid}). Scale columns of $B$ corresponding to $S$ by $2$ and call the new matrix $B'$, i.e. $B'=BD$ where $D$ is a diagonal matrix whose diagonal entries corresponding to $S$ are $2$ and are $1$ otherwise. Then $\ker(B')$ is $2$-regular because $B$ is totally unimodular. Note that $\icp(B')\leq 1$ because imbalanced circuits of $\ker(B')$ correspond to $\delta_{K_n}(T)$ such that $\delta_{\tilde{G}}(T)\cap S$ is a nonempty proper subset of $\delta_{\tilde{G}}(T)$ where $\emptyset\neq T\subset V(\tilde{G})$ (see \cite[Proposition 2.3.1]{Oxley2011Matroid} for a description of circuits of cographic matroids as minimal edge cuts). Any two such imbalanced circuits pairwise intersect because $\tilde{G}$ is complete. 

Let $b=Bc$. We claim that $\{x\in\mathbb{Z}^m:B'x=b, x\geq 0\}\neq \emptyset$ if and only if the parity potential problem $(G,c,S)$ is a yes-instance. Indeed, we have:
\begin{align*}
    &\exists x\in \mathbb{Z}^m \text{ such that } B'x=b,x\geq 0 \\
    \Longleftrightarrow&\exists x\in \mathbb{Z}^m \text{ such that } BDx=Bc,x\geq 0 \\
    \Longleftrightarrow & \exists x\in \mathbb{Z}^m \text{ such that } Dx-c\in\ker(B),x\geq 0 \\
    \Longleftrightarrow &\exists x\in \mathbb{Z}^m \text{ such that } Dx-c\in r.s.(A),x\geq 0 \\
    \Longleftrightarrow &\exists \pi:V(G)\to\bb{R} \text{ such that } \forall e=(v,u)\in E(G), c^\pi(e):=c(e)+\pi(v)-\pi(u)\in\bb{Z}_{\geq 0} \text{ and } \forall e\in S, c^\pi(e)\in 2\Z  \\
    \Longleftrightarrow &\exists \pi:V(G)\to\bb{Z} \text{ such that } \forall e=(v,u)\in E(G), c^\pi(e):=c(e)+\pi(v)-\pi(u)\in\bb{Z}_{\geq 0} \text{ and } \forall e\in S, c^\pi(e)\in 2\Z 
\end{align*}
The last step holds because if $\pi:V(G)\to\mathbb{R}$ satisfies $c^\pi(e)\in \mathbb{Z}$ for all $e=(v,u)\in E(G)$, then $\pi(v)-\pi(u)\in\mathbb{Z}$, then connectivity of $\tilde{G}$ implies that the $\pi$-values of any two vertices differ by an integer, so we can add a constant to $\pi$ to make $\pi$ integral. \qed

We observe the following failure of Erd\H{o}s--P\'osa property for cographic matrices, which can be compared with \Cref{EP3}:
\begin{corollary}\label{EP4}
    There exist matrices $\{A_n\}_{n\geq 4, n\in\mathbb{Z}}$ such that for all $n$, $A_n$ represents cographic matroids, $\ker(A_n)$ is $2$-regular, $\icp(A_n)=1$ and $\ict(A_n)=n-1$ .
\end{corollary}
\begin{proof}
    Consider $M(K_n)^*$ for $n\geq 4$. Circuits of $M(K_n)^*$ are minimal edge cuts of $K_n$, i.e. edge sets $\delta(V)$ for $\emptyset\neq V\subset V(K_n)$ (\cite[Proposition 2.3.1]{Oxley2011Matroid}). $M(K_n)^*$ are regular matroids, so there are TU matrices $M_n$ which represent $M(K_n)^*$. For each $n\geq 4$, pick a Hamiltonian path $P_n$ in $K_n$, scale columns of $M_n$ which correspond to edges on $P_n$ by $2$, and call these matrices $A_n$. An imbalanced circuit of $A_n$ is a circuit of $M(K_n)^*$ which intersects both $E(P_n)$ and $E(K_n)-E(P_n)$, which is an edge set $\delta(V)$ which intersects both $E(P_n)$ and $E(K_n)-E(P_n)$. But for all nonempty proper subset $V$ of $V(K_n)$, $\delta(V)$ intersects both $E(P_n)$ and $E(K_n)-E(P_n)$ because $(V(K_n), E(P_n))$ and $(V(K_n), E(K_n)-E(P_n))$ are connected graphs. So imbalanced circuits of $A_n$ are just $\delta(V)$ for $\emptyset\neq V\subset V(K_n)$. Any two such imbalanced circuits intersect, so $\icp(A_n)=1$. On the other hand, if $T\subseteq E(K_n)$ is an imbalanced circuit transversal, then $T\cap \delta(V)\neq\emptyset$ for all $\emptyset\neq V\subset V(K_n)$, so $(V(K_n), T)$ is connected, so an imbalanced circuit transversal of minimum size is a spanning tree, so $\ict(A_n)=n-1$.
\end{proof}

\section{Vertex Cover}
We mentioned in \Cref{intro} the connection between vertex cover and our matrix parameters $\icp$ and $\ict$. In this section we state these connections and discuss some open problems.

As in \Cref{intro}, let $G$ be an undirected graph, $A_G$ be the vertex-edge incidence matrix of $G$, and $M_G=[A_{G}^t| -I_{m}]$ where $m=|E(G)|$. $M_G$ is the constraint matrix appearing in the standard form IP formulation of the vertex cover problem. Below, $\oct(G)$ stands for \emph{odd cycle transversal number} of $G$, the minimum number of vertices which intersect every odd cycle of $G$. $\iocp(G)$ stands for the \emph{induced odd cycle packing number}, the maximum number of vertex-disjoint odd cycles of $G$, such that there is no edge with endpoints on two of these cycles.
\begin{lemma}\label{vc}
We have $\icp(M_G)=\iocp(G)$, $\ict(M_G)\leq \oct(G)\leq 2\cdot\ict(M_G)$, and $\ker(M_G) \text{ is 2-regular}$.
\end{lemma}
    For a proof, see \Cref{vcproof}. The main idea of the proof is to characterize the circuits of $M_G$, and connect imbalanced circuits of $M_G$ with odd cycles of $G$.
    
    While vertex cover is an NP-hard problem, it is fixed parameter tractable when parametrized by the size of a given odd cycle transversal. The running time has single exponential dependence on the parameter, for example $O^*(1.6181^k)$ for unit cost by \cite[Corollary 15]{NarayanaswamyEtAl2012LP}. \Cref{prop:milp} can be viewed as a generalization of FPT results for vertex cover parametrized by $\oct$ because $\ict(M_G)=\Theta(\oct(G))$ by \Cref{vc}, but the running time of \Cref{prop:milp} is $O^*(2^{O(k\log(k+1))})$, where the parameter dependence is slower than vertex cover parametrized by $\oct(G)$. While it can be very hard to improve general mixed-integer linear programming running time, there might be single exponential parameter dependence for \eqref{eq:ip} when $A$ (or $b$) is more general than having the form $M_G$ (or $\mathbf{1}_{E(G)}$). For example, since $\ker(M_G)$ is $2$-regular, we could ask whether \eqref{eq:ip} is solvable in time $O^*(2^k)$ if $\ker(A)$ is $2$-regular and we are given an imbalanced circuit transversal of size $k$.

    On the packing side, we have shown that \eqref{eq:ip} is fixed parameter tractable  parametrized by $\icp(A)$ when $A$ represents a graphic matroid. While \Cref{thm:cographic} shows cographic hardness, FPT algorithms parametrized by $\icp(A)$ might still exist for other problems, for example, vertex cover. It is  open whether vertex cover is NP-hard or polynomial time solvable on graphs with $\iocp(G)\leq 1$ (see introduction of \cite{BonamyEtAl2024}). In formulation \eqref{eq:ip}, this problem corresponds to $A=M_G, b=\mathbf{1}_{m}$, $c=(w,0)$ where $w\in\Q^{V(G)}_{\geq 0}$, and $\icp(A)\leq 1$. Positive progress on the class of the vertex cover problem with $\iocp(G)\leq 1$ could imply that $\icp$-parametrization of \eqref{eq:ip} works for a class of matrices larger than those representing graphic matroids.










\section*{Acknowledgments}
I am grateful to László Végh and Matthias Kaul for their mentorship, guidance, and support throughout the development of this paper.

\section*{Funding}\label{}
This research did not receive any specific grant from funding agencies in the public, commercial, or not-for-profit sectors.

\section*{Declaration of generative AI and AI-assisted technologies in the manuscript preparation process}
During the preparation of this work the author used ChatGPT Codex to assist with language editing, reference checking, and developing and clarifying some mathematical arguments. After using this tool/service, the author reviewed and edited the content as needed and takes full responsibility for the content of the published article.
\printcredits

\bibliographystyle{unsrtnat}

\bibliography{cas-refs}


\newpage
\appendix
\crefalias{section}{appendix} 
\section{Proof of \Cref{prop:milp}}\label{Prop3Proof}
    WLOG, suppose $S$ is the last $k$ columns of $A$ and let $T=[n]\backslash S$ be the first $n-k$ columns. By Gaussian elimination applied to $[A_T | A_S | b]$, we can transform $A_T$ to a totally unimodular matrix. Note that $A_S$ and $b$ might not be integral, as it can contain rational entries after Gaussian elimination. 

    Consider $\Lambda = \{z\in\bb{Z}^{k}:b-A_Sz\in\bb{Z}^m\}$. We claim that we can in polynomial time either decide $\Lambda =\emptyset$ or find invertible $M\in\mathbb{Z}^{k\times k}$ and $z_0\in\mathbb{Z}^k$ such that $\Lambda = z_0+M\mathbb{Z}^k$. Assuming this, the original problem is infeasible if $\Lambda = \emptyset$. Otherwise, it becomes 
    \begin{equation*}
        \begin{split}
            &\quad \min\{c_T^t x_T +c_S^t x_S:A_Tx_T+A_Sx_S=b, x_T\in\bb{Z}^{n-k}_{\geq 0}, x_S\in\bb{Z}^{k}_{\geq 0}\}
            \\ &=c_S^t z_0+\min\{ c_T^t x_T + c_S^t Mz : A_T x_T+ A_SMz =b-A_Sz_0,Mz\geq-z_0, x_T\in\bb{Z}^{n-k}_{\geq 0}, z\in\bb{Z}^k\} 
            \\ &= c_S^t z_0+\min\{ c_T^t x_T + c_S^t Mz : A_T x_T+ A_SMz =b-A_Sz_0,  Mz\geq-z_0, x_T\in\bb{R}^{n-k}_{\geq 0},  z\in\bb{Z}^k\}
        \end{split}
    \end{equation*}
    \noindent
    where in the first equality we applied our description of $\Lambda$ to $x_S$, and in the second equality we relax integrality constraints on $x_T$ because for any fixed $z\in\bb{Z}^k$, the set of values that $x_T$ can take forms an integral polyhedron (as $A_T$ is TU and the right hand side is integral). We now need to solve a mixed integer linear program with $k$ integer variables. It can be solved in deterministic time $2^{\ca{O}(k\log (k+1))}\size^{\ca{O}(1)}$, by a result of Kannan (see page 434 of \cite{Kannan1987Minkowski}).

    We now show how to compute $\Lambda$. First find suitably large $q\in\bb{N}$ such that $qb$ and $qA_S$ are integral. Compute the Smith normal form of $qA_S$, i.e. find nonsingular unimodular $U\in\bb{Z}^{m\times m}$ and $V\in\bb{Z}^{k\times k}$ such that $U(qA_S)V=D$ where the first $r$ ($r\leq\min(m,k)$) entries on the diagonal of $D$ are nonzero integers, and all other entries of $D$ are zero. Now $\Lambda=\{z\in\bb{Z}^k: qb-qA_S z\in q\bb{Z}^m\} = \{z\in\bb{Z}^k: U(qb)-DV^{-1}z\in q\bb{Z}^m\}=V\{z\in\bb{Z}^k: U(qb)-Dz\in q\bb{Z}^m\}$. Let $b'=U(qb)$. If $q\nmid b'_{i}$ for any $i\in[m]\backslash [r]$, then $\Lambda=\emptyset$. Otherwise, we want to find $z\in\bb{Z}^k$ such that $d_{i,i}z_i\equiv b'_i\mod q$ for $i\in [r]$ and $z_i\in\bb{Z}$ for $i\in[k]\backslash [r]$. The congruence equation $d_{i,i}z_i\equiv b'_i\mod q$ has a solution if and only if $\gcd(d_{i,i},q)\mid b'_i$ (otherwise $\Lambda=\emptyset$), in which case all integral solutions are $z_i'+l\cdot \frac{q}{\gcd(d_{i,i},q)}$ for $l\in\bb{Z}$, where $z_i'$ is one solution. Let $w_0=(z_1',...,z_r',0,...,0)\in\bb{Z}^k$ and let $P\in\bb{Z}^{k\times k}$ be a diagonal matrix where $P_{i,i}=\frac{q}{\gcd(d_{i,i},q)}$ for $i\leq r$ and $P_{i,i}=1$ for $i>r$. Then $\{z\in\bb{Z}^k: U(qb)-Dz\in q\bb{Z}^m\}=w_0+P\bb{Z}^k$, so $\Lambda = Vw_0+VP\bb{Z}^k$, which is the form we want. Smith normal form can be computed in polynomial time, and congruence equations $d_{i,i}z_i\equiv b_i'\mod q$ can be solved in polynomial time by the extended Euclidean algorithm, so we can compute $\Lambda$ in polynomial time. \qed
\section{Proof of \Cref{vc}}\label{vcproof}
    Suppose $(y,s)$ is an elementary vector of $\ker(M_G)$ where $y$ is indexed by the vertices and $s$ is indexed by the edges. Then for all $e=\{u,v\}\in E(G)$, $s_e=y_u+y_v$. Note that for any vector in $\ker(M_G)$, its $s$-components are uniquely determined by its $y$-components. $\supp(y)\neq \emptyset$, because otherwise $(y,s)=0$. 
    
    We claim that there is a subtree $T$ of $G$ such that $V(T)=\supp(y)$ and $E(T)\cap\supp(s)=\emptyset$. Otherwise, there exists a nonempty proper subset $U\subset \supp(y)$ such that $s_e\neq 0$ for all $e\in \delta(U)\cap \delta(\supp(y)-U)$. Then define $(y',s')$ by setting $y'_v=y_v$ for $v\in U$, $y'_v=0$ elsewhere, and $s'_e=y'_u+y'_v$ for all $e=\{u,v\}\in E(G)$, we get $(y',s')\in\ker(M_G)$ with strictly smaller support than $(y,s)$, a contradiction. 

    Take any $r\in V(T)$. Because $s(e)=0$ for all $e\in E(T)$, we see that $y(v)=y(r)$ for $v\in V(T)$ such that the $r$-$v$ path on $T$ has even length, and $y(v)=-y(r)$ for $v\in V(T)$ such that the $r$-$v$ path on $T$ has odd length.

    Conversely, we show that any vector in $\ker(M_G)$ defined by assigning its $y$-value alternately on a tree is an elementary vector. Let $T'$ be a subtree of $G$. pick any $r\in V(T')$ and let $y'_v=1$ if $v\in V(T')$ and the $v$-$r$ path on $T'$ has even length, let $y'_v=-1$ if $v\in V(T')$ and the $v$-$r$ path on $T'$ has odd length, and let $y'_v=0$ if $v\notin V(T')$. $s'$-values are then uniquely determined by $s'_{\{u,v\}}=y'_u+y'_v$. Then $(y',s')\in \ker(M_G)$. Furthermore, $(y',s')$ is an elementary vector of $\ker(M_G)$, because if $(y'',s'')\in\ker(M_G)\backslash\{0\}$ satisfies $\supp(y'',s'')\subseteq \supp(y',s')$, then $y''(v)\neq 0$ for some $v\in \supp(y')=V(T')$. Then since $s''(e)=s'(e)=0$ for all $e\in E(T')$, the values of $y''$ on $V(T')$ are uniquely determined up to scaling, so $(y'',s'')$ is a scalar multiple of $(y',s')$.

    Above all, any elementary vector $(y,s)$ of $\ker(M_G)$ can be described (up to multiplication by a nonzero scalar) as $(y^T,s^T)$ where $T$ is a subtree of $G$ and $(y^T,s^T)$ is constructed as in the previous paragraph. Note that $(y^T,s^T)$ can be scaled to a  $\{0, \pm 1, \pm 2\}$-vector, so $\ker (M_G)$ is $2$-regular.

    Moreover, if $C$ is an odd cycle of $G$, let $T$ be the path we get from removing one edge $e$ of $C$, then $(y^T,s^T)$ is an imbalanced elementary vector because $|s^T(e)|=2|y^T(v)|$. for all $v\in V(C)$. Then it is not hard to see that $\icp(M_G)\geq \iocp(G)$. For the other direction, we note that if $(y^T,s^T)$ is an imbalanced elementary vector, then $|s^T(e)|=2|y^T(v)|$ for some $e=\{u,w\}$ and $v\in V(T)$ where $u,w\in V(T)$ and the $u$-$w$ path on $T$ has even length. This path together with $e$ forms an odd cycle. Then it is not hard to see that $\icp(M_G)\leq \iocp(G)$.

    We have shown that for any imbalanced elementary vector $(y^T,s^T)$, $V(T)$ contains the vertex set of an odd cycle, so $\ict(M_G)\leq \oct(G)$. 

    If $A\dot\cup B $ is an imbalanced circuit transversal of $M_G$ where $A\subseteq V(G)$ and $B\subseteq E(G)$, let $S=A\cup\bigcup_{e=\{u,v\}\in B}\{u,v\}$, then $|S|\leq 2|A\dot\cup B|$ and $S$ is an odd cycle transversal of $G$ because from any odd cycle of $G$ we can build an imbalanced elementary vector $(y^T,s^T)$ where $T$ is the odd cycle without one edge. So $\oct(G)\leq 2\cdot\ict(M_G)$. \qed

\end{document}